*A model-based technique to identify lubrication condition of hydrodynamic bearings using the rotor vibrational response*

Marcus Vinícius Medeiros Oliveira*

Bárbara Zaparoli Cunha

Gregory Bregion Daniel

Corresponding author email: m263128@dac.unicamp.br

School of Mechanical Engineering, University of Campinas, Brazil

*Abstract*

Faults related to hydrodynamic bearing can imply in high maintenance costs when late-detected and even to the total shutdown of the system. Thus, techniques of early fault diagnosis have high relevance to the reliability of rotating machinery. However, a common fault caused by inadequate bearing oil supply has not yet received appropriate attention. This paper presents a new approach to model and identify starved or excessive oil supply in hydrodynamic bearings. The developed identification technique is a model-based process that uses the rotor vibration signal to access the bearing lubrication. Numerical identification tests were performed and the results showed that the proposed method can satisfactorily estimate the oil flow rate in bearings under starved and flooded lubrication conditions, thus representing a useful and promising tool for condition monitoring and fault diagnosis applied to rotating machinery.

*1. Introduction*

Hydrodynamic journal bearings have wide applicability in rotating machinery, being a crucial component to ensure good performance and safety of the operation. For this reason, it is mandatory predicting the influence of bearings conditions on rotors' dynamic behavior in order to develop more reliable and efficient machinery. In addition to contributing to the design phase, understanding bearing dynamic phenomena and how they affect the response of the rotating system is important for machinery health monitoring. The knowledge of vibrations signatures related to bearing conditions, commonly observed in the frequency spectrum, can indicate the occurrence of a fault in the component. Thus, early fault diagnosis can provide more effective maintenance programs, reducing downtime and preventing sudden breakdown, which can cause huge production losses, high repair costs and even environmental threats.

The fault diagnosis begins by detecting the fault occurrence and then goes to fault identification, which determines the type, magnitude and location of the fault [1]. In rotordynamics, model-based diagnosis and identification are a widely used technique, in which the fault conditions are represented by models from the parameters to be adjusted seeking to obtain a simulated response, e.g. lateral vibration, similar to that measured on the real machine

In this context, several works have been addressed to hydrodynamic bearing fault identification using model-based methods. Bachschmid et al. [2] presented a model-based identification method for multiple rotor faults modeled as forces and moments acting on the rotating system, being one of them the ovalization of the journal bearing. The identification is carried out by least-squares fitting, minimizing the multidimensional residual between the vibrations measured in the machine and the vibrations calculated with the faults. Recently, Alves et al. [3] proposed to identify the bearing ovalization by combining model and machine learning strategies. Using the output data from an ovalized bearing model, a convolutional neural network (CNN) algorithm was trained to identify bearing ovalization fault conditions under different operational conditions, resulting in a good fault diagnosis process.

Diagnosis of wear on hydrodynamic bearings has also been extensively investigated. Papadopoulos et al. [4] proposed a method to identify worn bearing radial clearance by minimizing the sum of squared differences between the measured and simulated shaft responses. Gertzos et al. [5] proposed a graphical method for bearing wear depth identification based on the intersection of diagrams of bearing dynamics characteristics measured or calculated from the Somerfield number and wear depth. Machado and Cavalca [6] identified wear depth and angular position in journal bearings by minimizing a logarithmic least-squares function that compares the measured and the simulated rotor's unbalance frequency response in directional coordinates. In continuity to that work, Mendes et al. [7] used the same method to identify bearing wear, but comparing the directional frequency response functions (dFRF) instead of the directional unbalance frequency responses. Later, Alves et al. [8] performed bearing wear identification comparing only the vibration information from shaft's displacement inside the bearing, in which the measured and simulated vibrations responses were compared in terms of forward and backward components in the frequency domain.

Another critical bearing fault observed in rotating machinery is related to the oil starvation. Despite being a common fault condition in hydrodynamic bearings, it has not yet been the object of model-based identification strategies. The starvation condition may occur when the oil supply is insufficient to fill the radial clearance between the bearing and the shaft surfaces in the groove region [9]. This fault can occur due to wrong setting of the oil flow rate, problems with the oil pump and/or measurement sensors (flow and/or pressure) and also due to leaks or clogging in the oil supply line. A starved hydrodynamic bearing has a modified pressure field in the oil film and, consequently, a different performance. Artiles and Heshmat [10] showed that smaller oil flow rates reduce the extent of the hydrodynamic pressure region, affecting rigidity and damping of the oil film. Vincent et al. [11] verified a drastic reduction in load capacity of a bearing under severe starvation condition. Tanaka [12] showed that the shaft eccentricity can be severely increased under starvation condition. On the other hand, flooded lubrication with excessive oil flow rates may worsen the bearing performance and, in some cases, it may be desirable to operate with lower flowrates. Vijayaraghavan et. al [13] showed that oil supply can be reduced in order to save power loss of power, but it may not be safe operating with supply rates below 70 per cent of flooded condition. Hashimoto and Ochiai [14] suggested that starved lubrication can be used to stabilize a small-bore journal bearing. In this context, fundamental studies are still being done to better understand conformal contacts under starvation and optimize lubricant supply in hydrodynamic bearings [15,16].

The studies previously mentioned have shown the importance of oil supply on hydrodynamic bearings behavior, and consequently, on the rotating system performance. For this reason, monitoring systems able to identify the oil supply conditions in hydrodynamic bearings can provide important information to ensure good performance and operational safety for the rotating machinery. Recently, Poddar and Tandon [17] analyzed starvation in journal bearings by measuring acoustic emission (AE) and vibration. Later, the same authors proposed the use of AE in lubrication conditions monitoring through a machine learning approach [18]. Another way to monitor oil flow conditions that has not yet been explored is to perform identification using model-based methods. Based on this, the proposed work aims at using the rotor's vibration and a model-based method to identify the oil flowrates, allowing the detection of starvation conditions or even excessive oil supply. For this, a mass-conserving bearing model has been applied considering the oil flowrate as an input variable. Thus, the identification is performed by adjusting the oil flowrate in the model of the hydrodynamic bearings in order to minimize the error between the simulated and reference responses. The simulated vibrational responses of the rotor are analyzed in terms of backward and forward components in the frequency domain, using the orbit angle and the forward-backward ratio as response-parameter to the identification process. The analyses performed in this work show that the proposed approach is able to successfully identify the oil flowrate of the hydrodynamic bearings from the dynamic response of the rotor, indicating conditions of faults by starvation. Therefore, the development of a new technique for oil flowrate identification is the main contribution of this work, representing a useful and promising tool for condition monitoring and fault diagnosis applied to rotating machinery.

## 2. Methodology

This chapter presents the formulation related to the proposed approach in this work. First, Section 2.1 describes the hydrodynamic bearing model. In following, Section 2.2 presents the model of the rotating system. Lastly, Section 2.3 describes the numerical approach used to the oil supply flowrate identification.

### 2.1 Bearing Model

In order to simulate the influence of the bearing oil supply condition on the rotor dynamic response, it is necessary to have a bearing model that considers the inlet oil flow in the bearing groove. Figure 1 shows a schematic representation of hydrodynamic journal bearing with oil groove, in which R is the shaft radius, L is the bearing width, a is the groove width, b is the groove circumferential length, $e_X$ and $e_Y$ are the horizontal and vertical component of the shaft's eccentricity, α is the angular coordinate, Ω is shaft rotational speed and $h$ is the oil film thickness, as described in Eq. (1):

$$h(\alpha) = c_r + e_X \cdot \sin\alpha - e_Y \cdot \cos\alpha \quad (1)$$

where $c_r$ is the bearing radial clearance.

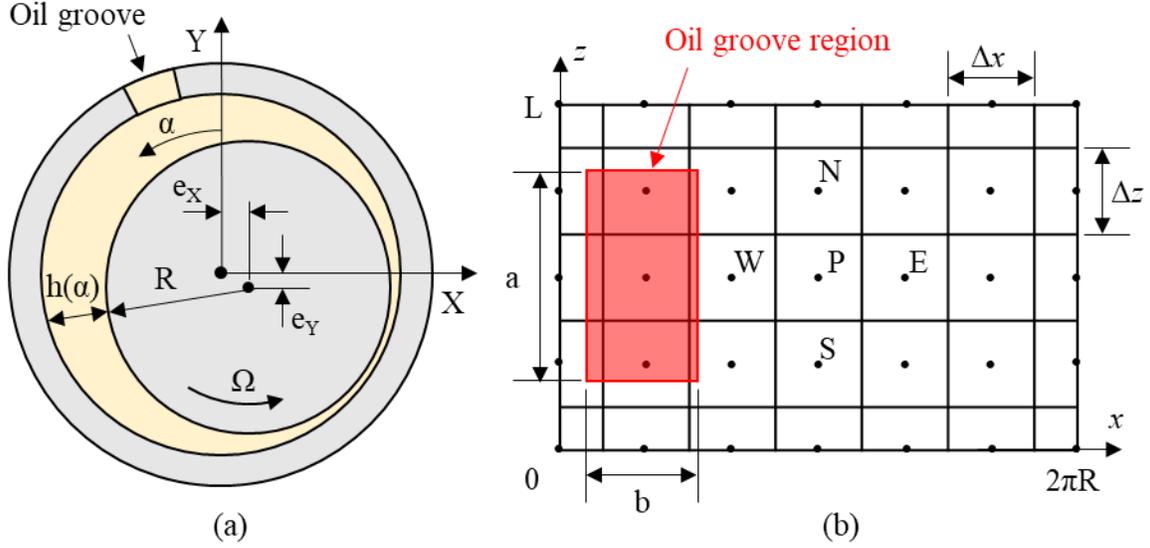

Fig.1. (a) Bearing parameters and (b) discretized oil domain.

Generally, the pressure distribution is modeled by classical Reynolds equation and the Swift-Stieber conditions to represent cavitation effects. However, the imposition of these conditions does not ensure the conservation of the mass of lubricant, making this approach unsuitable to consider the oil supply flow rate in the model. Therefore, this work proposes a new approach based on the solution by Finite Volume Method (FVM) of a mass-conserving lubrication model, also known as $p - \theta$ model [19]. In the $p - \theta$ model, the Reynolds equation is modified by introducing a scalar field of the oil fraction, denoted by $\theta$, as indicated in Eq. (2):

$$\frac{\partial}{\partial x}\left(\frac{h^3}{12\mu}\frac{\partial p}{\partial x}\right) + \frac{\partial}{\partial z}\left(\frac{h^3}{12\mu}\frac{\partial p}{\partial z}\right) = \frac{\partial}{\partial x}\left(\frac{U}{2}h\theta\right) + \frac{\partial(h\theta)}{\partial t} \qquad (2)$$

where $p$ is the pressure field, $x$ and $z$ are the circumferential and axial coordinates of the oil film, $t$ is the time, $\mu$ is the oil dynamic viscosity and $U$ is the linear speed on the shaft surface. According to Fig. 1(b), the oil film is discretized in the $x$ and $z$ directions and the FVM is applied to numerically solve the Eq. 2, thus obtaining the pressure and oil fraction distributions on the bearing. Using the central difference scheme for the evaluation of the partial derivatives in the Poiseuille terms (left side terms), the upwind scheme in the Couette terms (first right side term) and an explicit evaluation in the Squeeze term (second right side term), a system of equations $p - \theta$ for every volume $P$ in the mesh is obtained as:

$$\theta_P C_1 + p_P(C_3 + C_4 + 2C_5) = \theta_W C_2 + p_E C_3 + p_W C_4 + p_N C_5 + p_S C_5 \qquad (3)$$

in which $C_1$ to $C_5$ are coefficients from the integration of Eq. 2 and the subscripts E, W, N and S refer to adjacent volumes around the volume $P$ in directions east, west, north and south, respectively.

In order to be able to model the bearing lubrication as a function of the oil supply flow, a new approach is proposed considering the inlet flow as a boundary condition of the lubrication problem. Since the terms in Eq. 3 are evaluated as flux balance of lubricant over each finite volume, this can be done by adding a source term to Eq. (3) for the volumes in the oil groove region:

$$\theta_P C_1 + p_P(C_3 + C_4 + 2C_5) = \theta_W C_2 + p_E C_3 + p_W C_4 + p_N C_5 + p_S C_5 + \frac{Q_s}{N_g} \qquad (4)$$

being $Q_s$ the total supply flow rate and $N_g$ the number of volumes inside the groove.

To obtain the pressure and the oil fraction distributions, the system of equation is solved iteratively for the entire mesh, applying the Gauss-Seidel method with the boundary conditions:

$$\begin{cases} p_{(z,0)} = p_{(z,2\pi R)} \\ \theta_{(z,0)} = \theta_{(z,2\pi R)} \\ p_{(0,x)} = p_{(L,x)} = 0 \\ \theta_{(0,x)} = \theta_{(L,x)} = 1 \end{cases} \qquad (5)$$

During the iterative procedure, the values of $p$ and $\theta$ must respect the cavitation condition in every volume:

$$\begin{cases} p_P > 0 \rightarrow \theta_P = 1 \\ \theta_P < 1 \rightarrow p_P = 0 \end{cases} \qquad (6)$$

Then, the converged pressure field is integrated along the bearing domain to obtain the hydrodynamic forces. However, the simulation of the rotor dynamic response would require the recalculation of the pressure field for each time step what would lead to high computational cost. In order to reduce the computational cost and knowing that the shaft's displacements are small, the hydrodynamics forces can be linearized around the equilibrium position, as proposed by Lund [20]. Thus, the hydrodynamic forces are representing by stiffness and damping coefficients, that are the partial derivatives in relation to the shaft's displacement and velocity evaluated at the equilibrium position $e_{X0}$ and $e_{Y0}$, as presented in Eq. (7):

$$F_{HX} = K_{XX}(e_X - e_{X0}) + K_{XY}(e_Y - e_{Y0}) + C_{XX}\dot{e}_X + C_{XY}\dot{e}_Y \qquad (7)$$

$$F_{HY} = K_{YY}(e_Y - e_{Y0}) + K_{YX}(e_X - e_{X0}) + C_{YY}\dot{e}_Y + C_{YX}\dot{e}_X$$

being $K$ the stiffness coefficient, $C$ the damping coefficient, $F_H$ the hydrodynamic forces, $e$ and $\dot{e}$ the shaft's eccentricity and translational velocity, $X$ and $Y$ the global coordinates on the inertial reference system, as presented in Fig. 1.

2.2 Rotating System Model

In this work, the rotating system is modeled using Finite Element Method (FEM), being composed of cylindrical Timoshenko's beam elements, rigid disc elements and two hydrodynamic bearings. This model has four degrees-of-freedom (DOF) per node, being two translational displacements and two angular displacements (lateral vibration), as proposed by Nelson and

McVaugh [21]. Assembling the equation of motion of all elements that compose the rotating system, the global equation of motion can be written as presented in Eq. (8).

$$[M]\{\ddot{q}(t)\} + ([C] + \Omega[G])\{\dot{q}(t)\} + [K]\{q(t)\} = \{F(t)\} \tag{8}$$

where $[M]$, $[C]$, $[G]$ and $[K]$ are the global mass, damping, gyroscopic and stiffness matrices, respectively; $\{q(t)\}$ is the displacement vector of the node's degrees of freedom and $\{F(t)\}$ is the acting forces vector. For a rotating system, the damping matrix can be considered proportional to the stiffness matrix. The forces vector includes gravitational force (shaft's weight), bearing hydrodynamic forces ($F_{H_X}$ and $F_{H_Y}$), and unbalance force $F_D$, that is modelled as presented in Eq. (9):

$$\{F_D(t)\} = \begin{Bmatrix} m_u \Omega^2 \cos(\Omega t) \\ m_u \Omega^2 \sin(\Omega t) \\ 0 \\ 0 \end{Bmatrix} \tag{9}$$

where $m_u$ is the unbalance moment.

The dynamic response of the rotating system can be obtained by solving the equation of motion presented in Eq. (8). Thus, the shaft's orbit is composed of vibrations in X and Y directions, that can be represented in the complex plane as:

$$Z(t) = V(t) + iW(t) \tag{10}$$

Being $i = \sqrt{-1}$, and $V$ and $W$ are the displacements in X and Y directions, respectively.

In the complex plane, the vibration signal can be represented as a complex harmonic signal resultant of contra-rotating forward and backward phasors, also called directional response [22]. Applying the Discrete Fourier Transform to the complex vibration signal, the amplitude and phase of the forward and backward components are distinguished in the frequency full spectrum, which aids the recognition of fault signatures on rotor dynamic behaviour [8]. For this reason, the amplitude and phase of the rotor directional response are used in the identification method proposed in this work.

As the rotating system model is linear, the simulated response is predominantly composed of synchronous components. Thus, the amplitudes and phases are evaluated at the positive and negative shaft's rotational frequency in the full spectrum, being -1X for backward component and 1X for forward component.

*2.3 Oil Supply Flowrate Identification*

In the previous section, the rotating system model is presented, being composed of a bearing model that considers the oil supply flowrates as input parameters. Based on this model, this work proposes an approach for oil supply flowrate identification. For this, the flowrates are adjusted in the model in order to obtain a simulated response similar to a measured response in the rotating system whose oil flows are unknown. For this purpose, an identification process is proposed considering a root-finding problem, in which the equation represents the relative error between the measured and simulated responses of the rotating system. Thus, the proposed identification approach aims to find the combination of oil flowrates that is the root of a scalar error function defined as:

$$\varepsilon(Q_1, Q_2) = \frac{f_s(Q_1,Q_2) - f_m(Q_1,Q_2)}{f_m(Q_1,Q_2)} \tag{11}$$

where $f_s$ and $f_m$ are scalar response-parameter derived from the simulated and measured directional responses, respectively, which depends on the oil supply flowrate on the bearing 1 ($Q_1$) and on the bearing 2 ($Q_2$).

In order to ensure a unique solution, the vibration information is evaluated at two points of the rotor, more specifically, at the two bearings. Therefore, the identification process becomes a search for the solution of the system of equations (error functions):

$$\begin{cases} \varepsilon_1(Q_1, Q_2) = \frac{f_{s_1}(Q_1,Q_2) - f_{m_1}(Q_1,Q_2)}{f_{m_1}(Q_1,Q_2)} = 0 \\ \varepsilon_2(Q_1, Q_2) = \frac{f_{s_2}(Q_1,Q_2) - f_{m_2}(Q_1,Q_2)}{f_{m_2}(Q_1,Q_2)} = 0 \end{cases} \tag{12}$$

where the subscripts 1 and 2 refer to the error function and its response-parameter evaluated at bearing 1 and 2, respectively.

As previously described, the directional response is analyzed in terms of amplitudes and phases of the forward and backward components. Hence, these components can be combined to generated a unique response-parameter for each bearing or they can be independently used, creating $n \geq 2$ error functions to compose the system of equations. In following, the system of equations is solved using the generalization of the Newton-Raphson method, in which the error functions are linearized and the roots $\{Q\} = [Q_1 \; Q_2]^T$ are iteratively obtained from an initial guess, $\{Q\}_{k=0}$, by computing:

$$\{Q\}_{k+1} = \{Q\}_k + \{s\}_k \tag{13}$$

being $\{s\}_k$ the solution of the linear system written in Eq. (14):

$$[J]_k \{s\}_k = -\{S\} = -\begin{Bmatrix} \varepsilon_1 \\ \vdots \\ \varepsilon_n \end{Bmatrix}_k \tag{14}$$

in which $\{S\}$ is the error functions and $[J]_k$ is the Jacobian of $\{S\}$.

The iterative process stops when the sum of the absolute value of all error functions ($\varepsilon_T$), is sufficiently small. In order to improve the robustness of the proposed method, the search domain is constrained considering an interval of reasonable flowrate for each bearing. Furthermore, a multi-start scheme with 9 initial guesses is used from a 3x3 grid and the number of iterations $k$ is limited per initial guess in the search process. Then, $\varepsilon_T$ is evaluated at the 9 initial guesses and the identification (search process) starts from the guess with the first lowest error. If the convergence is not reached, the process restarts from the guess with second lowest error, and so on. The identification process is presented in Fig. 2.

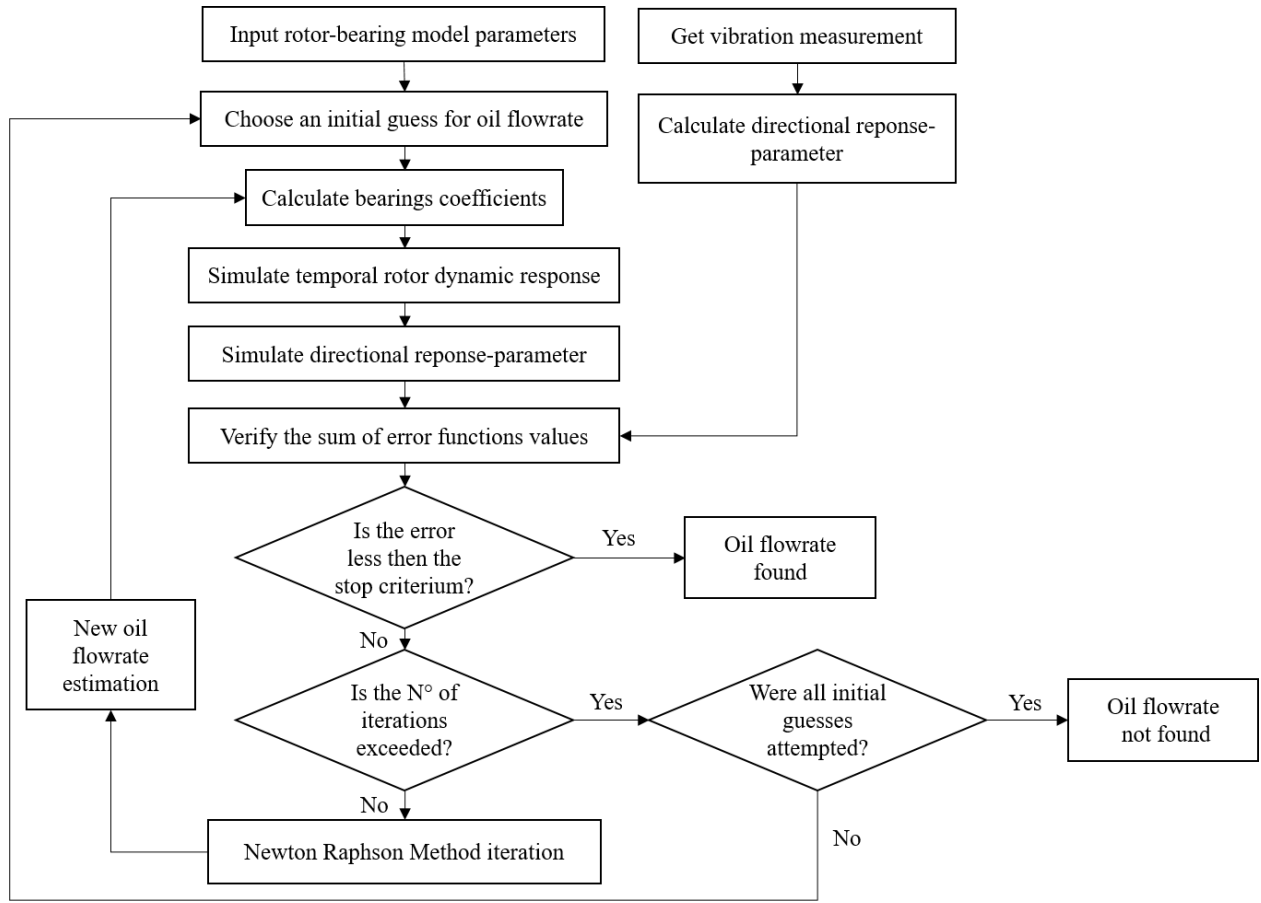

Fig. 2. Schematic flowchart of the identification process.

*3. Numerical Results*

This section presents the numerical results of oil flowrate identification for different lubrication conditions. The rotor-bearings model used in this work is shown in Fig. 3. The rotating system is composed of a shaft discretized into 20 beam elements (21 nodes), three disc elements (at nodes 10, 13 and 15) and two identical cylindrical bearings (at nodes 6 and 20). The shaft elements' dimensions are described in Table 1, while the disc elements' dimensions are presented in Table 2. The material of the discs and shaft is steel with Young Modulus E = 210GPa, Poisson ratio v=0.3 and density $\rho = 7850$ kg/m³. The only source of excitation is an unbalance moment $m_u = 1.24 \times 10^{-2}$ Nm acting on the central disc. The total rotor weight $W_y = 6545.6$N is equally supported by the two bearings whose characteristics are presented in Table 3. The bearings were modeled with a finite volume mesh composed of 100 volumes in both circumferential and axial direction, being chosen after mesh convergence tests.

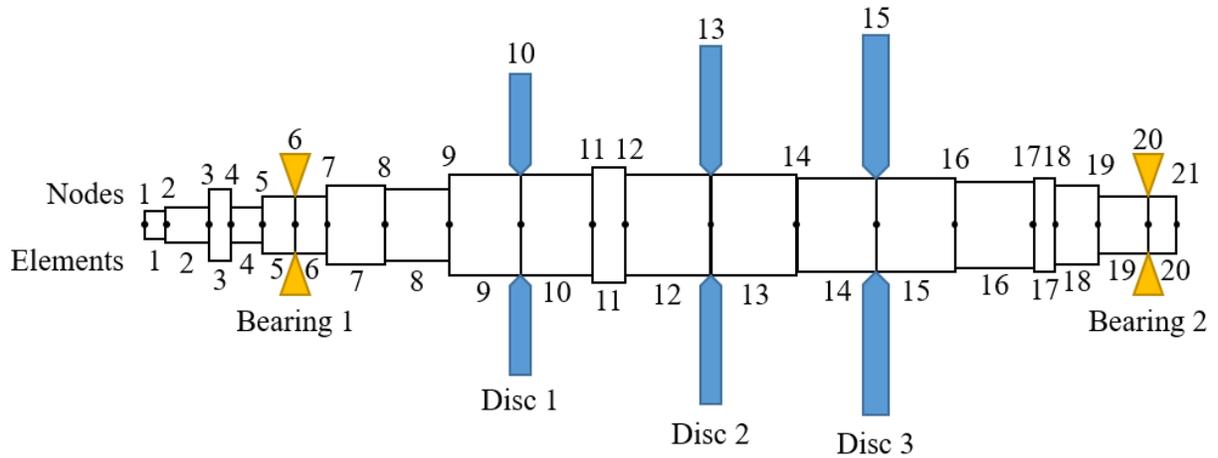

Fig. 3. The finite element model of the rotor-bearing system.

Table 1

Shaft element dimensions.

| Element Number | Length [mm] | Diameter [mm] |
| --- | --- | --- |
| 1 | 30.0 | 30.0 |
| 2 | 60.0 | 52.5 |
| 3 | 30.0 | 112.5 |
| 4 | 45.0 | 52.5 |
| 5,6 | 45.0 | 90.0 |
| 7 | 120.0 | 135.0 |
| 8 | 135.0 | 112.5 |
| 9,10 | 180.0 | 187.5 |
| 11 | 45.0 | 262.5 |
| 12,13 | 240.0 | 187.5 |
| 14,15 | 210.0 | 165.0 |
| 16 | 210.0 | 150.0 |
| 17 | 30.0 | 165.0 |
| 18 | 60.0 | 135.0 |
| 19 | 78.4 | 90.0 |
| 20 | 41.6 | 90.0 |

Table 2

Disc element dimensions.

| Disc | Width [mm] | External Diameter [mm] |
| --- | --- | --- |
| 1 | 50.0 | 525.0 |
| 2 | 50.0 | 600.0 |
| 3 | 50.0 | 697.5 |

Table 3

Bearing and lubricant oil properties.

| Parameters and units | Value |
|---|---|
| Radius [mm] | 45.0 |
| Width [mm] | 70.0 |
| Radial clearance [µm] | 120 |
| Oil dynamic viscosity [Pa.s] | 0.094 |
| Groove angular position [°] | 0 |
| Groove circumferential length [mm] | 16.2 |
| Groove width [mm] | 35.0 |

Firstly, the sensitivity of the directional response in relation to oil supply flowrate variation was evaluated. Then, numerical identifications for different flowrates were carried out, encompassing flooded and starved lubrication conditions. The performance of the identification method is evaluated based on the errors obtained and also the processing time. The analyses were conducted at a rotating frequency of 75 Hz (operation rotational speed), which lies between the first critical speed and the instability threshold of the system. The rotor dynamic response was simulated for a time of 15 s with time step of $1 \times 10^{-4}$ s. The first 10 s of the simulated response was discarded to avoid numerical transient and the last 5 s was used for the identification process of the oil supply flowrate on the bearings.

*3.1 Sensitivity of directional response*

The proposed identification method requires that the system presents distinguishable responses for different oil supply flowrates. Thus, the sensitivity of the directional components in relation to the flowrate variation was verified from 50% to 150% of the threshold for flooded lubrication ($Q_T = 546$ ml/min), i.e. for supply flowrates varying from 273 to 819 ml/min.

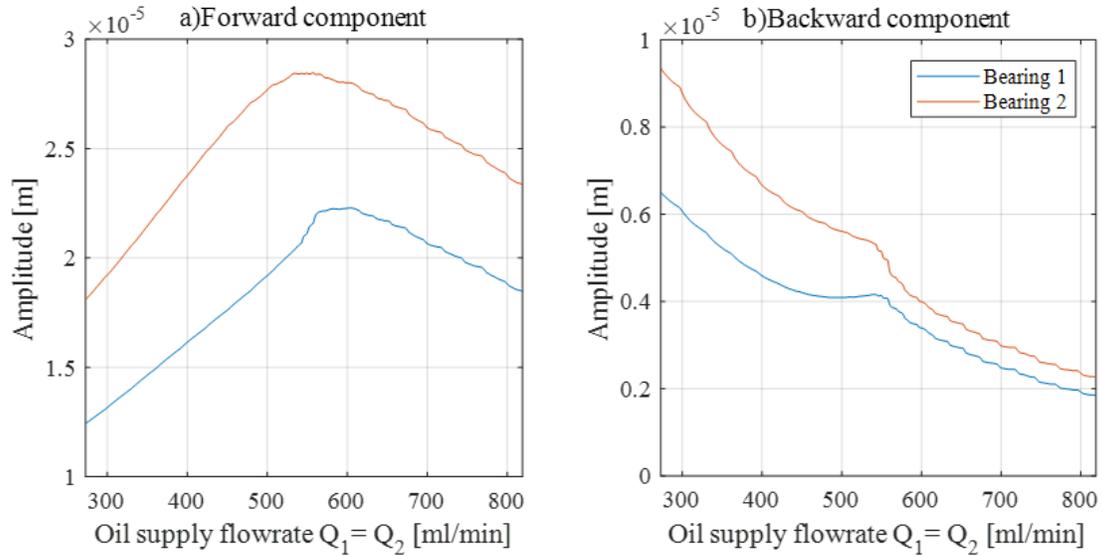

Fig. 4. a) Forward and b) backward amplitudes variation.

Fig. 4 shows the influence of the oil supply flowrate on the amplitudes of the forward 1X and backward -1X components from the shaft's displacement inside the bearings 1 and 2 (nodes 6 and 20, respectively). Since the rotor geometry is not symmetric, the dynamic responses are different inside the bearings 1 and 2. However, the behavior of these responses are similar and clearly sensitives to the flowrate variation. Observing Fig. 4a, the forward amplitude variation presents a critical point around $Q_T$. Thus, if only the forward amplitude is used as a response-parameter to create an error function, this function would have two roots because Eq.11 is a vertical translation of the response-parameter curve. From Fig. 4b, the backward amplitude variation is monotonic at bearing 2, but presents a plateau near $Q_T$ for the response from bearing 1, also implying in error functions with multiple roots in this region. Obviously, multiple roots are not desirable in a root-finding problem, hence using forward or backward amplitudes apart may not be a good option for the identification of bearing oil supply flowrate. On the other hand, the variation of the forward and backward amplitudes is smooth, keeping the trend of the curve even for relatively small flowrate variation (about 1% of $Q_T$), especially for flowrates under starvation condition. Thus, error functions using these amplitudes as response-parameters present roots defined with enough precision to detect any practical changes in the bearing conditions. In this context, error functions that jointly use the forward and backward amplitudes as response-parameters tend to represent a promising alternative.

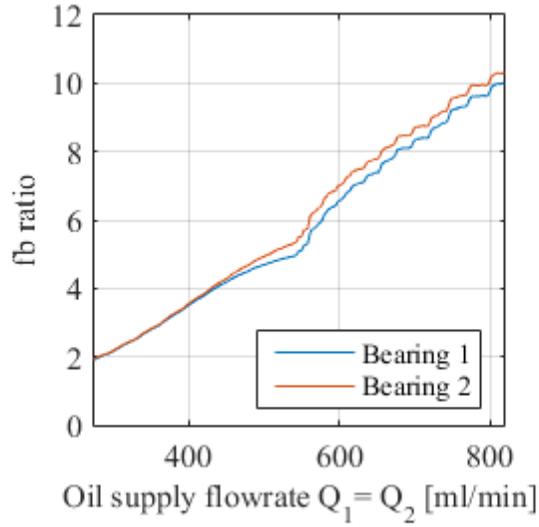

Fig. 5. Variation of the *fb* ratio.

Based on the observations discussed previously, different combinations of forward and backward signals were tested to take advantage of their sensitivity and smoothness. A suitable combination was obtained by the ratio between forward and backward amplitude, denoted *fb*, shown in Fig. 5. The *fb* curve has a monotone behavior, which eliminates the possibility of having an error function with multiple roots. Moreover, the *fb* ratio is quite linear and keeps the smoothness of the amplitudes, especially below $Q_T$.

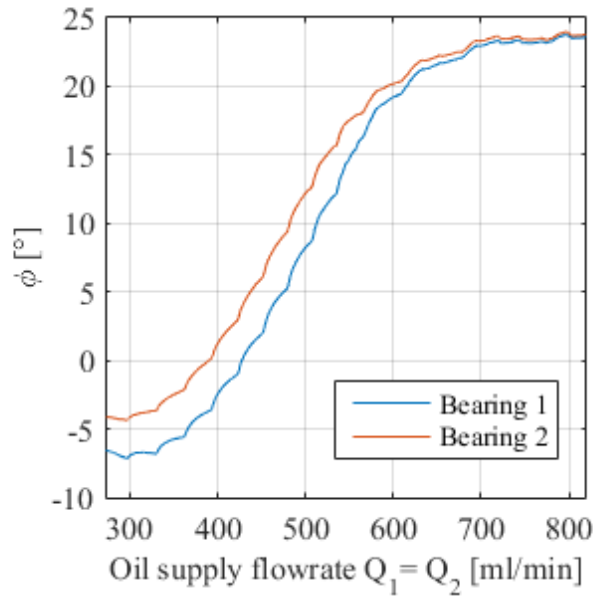

Fig. 6. Variation of the orbit angle $\varphi$.

Another response-parameter explored was the angle of the orbit $\varphi$ that is calculated from the arithmetic mean of the forward and backward phases. This parameter is also sensitive to the flowrate variation, as shown in Fig. 6. Over most of the flowrate range, $\varphi$ is monotonic and smooth enough to be used in the proposed flowrate identification. However, below 350 and above 700 ml/min, the angle $\varphi$ becomes nearly constant and its error function would present several roots.

As previously presented, $fb$ and $\varphi$ show to be suitable parameters to identify the oil flowrate when it is equal in both bearings. In order to verify if these response-parameters are still suitable to identify distinct oil flowrates in the bearings, it is necessary to observe the parameters' behavior when the oil supply flowrates vary independently in bearing 1, $Q_1$, and bearing 2, $Q_2$. Thus, Fig. 7 shows the contour plot of $fb$ and $\varphi$ as functions of $Q_1$ and $Q_2$, in which $fb_1$ and $\varphi_1$ are evaluated in the center of bearing 1 and $fb_2$ and $\varphi_2$ are evaluated in the center of bearing 2. As expected, the parameters are more sensitives to the variation of the flowrate in the bearing where they are evaluated. However, there is a strong dependency between the two bearings that makes the $fb$ and $\varphi$ sensitives to $Q_1$ and $Q_2$. This dependency tends to be reduced in more flexible rotors.

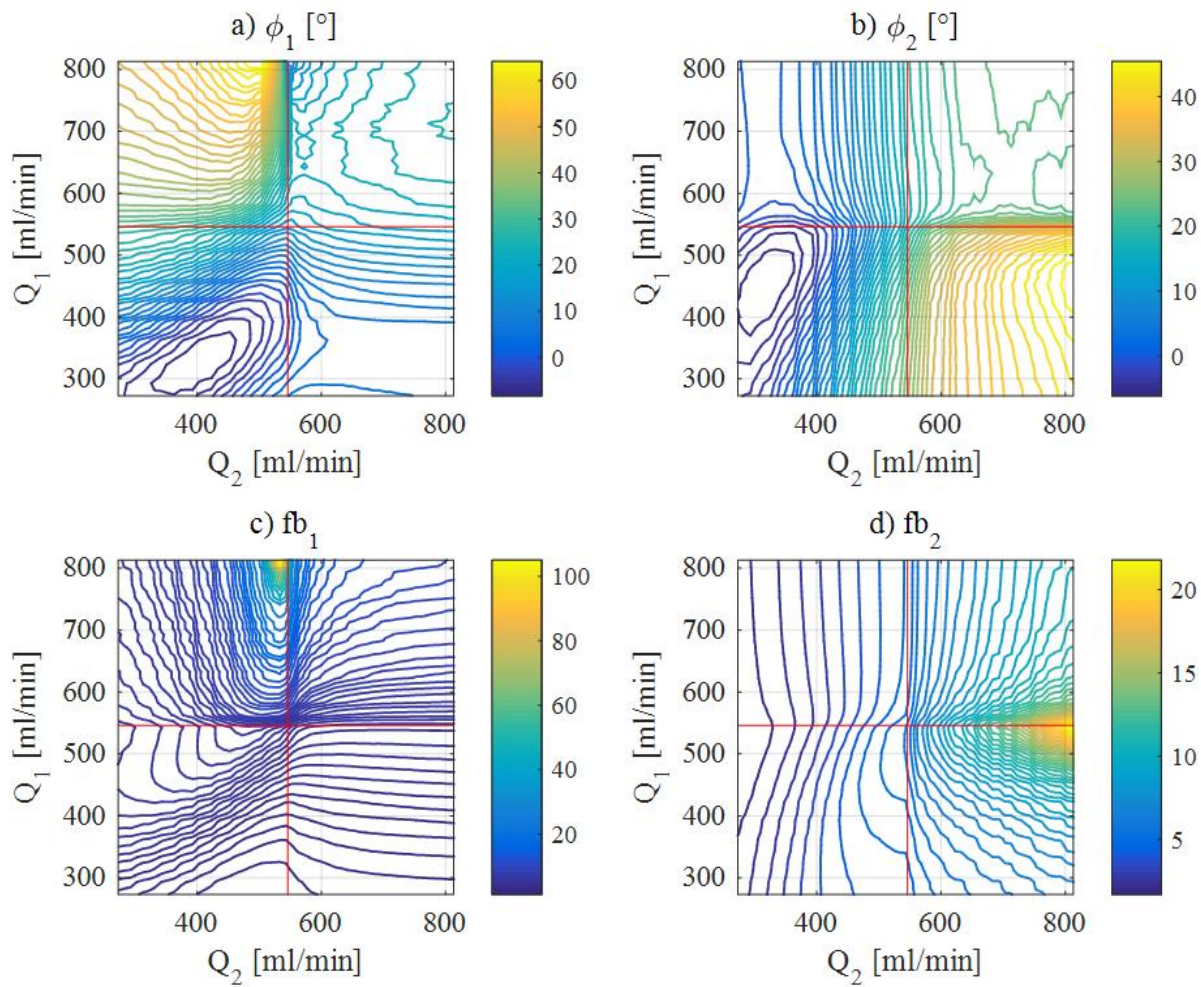

Fig. 7. Variation of $fb$ and $\varphi$ in bearings 1 and 2. Red lines indicate the threshold flowrate.

Graphically, unknown flowrates $Q_1$ and $Q_2$ can be identified at the intersection of *fb* or $\varphi$ levels observed in the system. Based on the contours presented in Fig. 7, it was observed that using *fb* or $\varphi$ from different bearings facilitates the intersection identification since the levels of the parameters from one bearing tend to be perpendicular to the levels of the parameters from the other bearing. In addition, *fb* shows to be suitable for identification because its contour plots have only one line per level, resulting in a unique solution. On the other hand, contour plots of $\varphi$ show closed lines per level, as observed in the third quadrant in Figs. 7a and 7b, which may result in multiple solutions. Thus, using error functions of $fb_1$ and $fb_2$ seems to be preferable for flowrates identification. However, the parameter $\varphi$ is useful as complementary information, for example, in cases near $Q_T$ where the $fb_1$ and $fb_2$ levels tend to be parallel. For this reason, the error functions proposed in this work are composed of both response-parameters, namely, the ratio between forward and backward amplitudes (*fb*) and angle of the orbit ($\varphi$).

## 3.2 Oil supply flowrate identification

From the results observed in the last section, both response-parameters *fb* and $\varphi$ were chosen to be used in the identification process. Thus, the system of error equations (Eq. 12) to be solved is composed of four error functions: $\varepsilon_{fb1}$, $\varepsilon_{fb2}$, $\varepsilon_{\varphi 1}$ and $\varepsilon_{\varphi 2}$. In order to test the proposed identification method under flooded and starved conditions, 8 reference oil flowrates between 75 and 125% of $Q_T$ (409.5, 448.5, 487.5, 526.5, 565.5, 604.5, 643.5 and 682.5 ml/min) were chosen to be identified. Thus, combining these reference oil flowrates in the bearings 1 and 2, a total of 64 cases were evaluated. The experimental response in each case was replaced by a reference response simulated with the known reference flowrates in the bearing 1 and 2. To better represent a real system, the reference response was contaminated with white-noise created with the MATLAB function *awgn* considering a Signal-to-Noise Ratio (SNR) of 12 dB, in agreement with the literature [3]. The noise was added to shaft's displacement, before the extraction of forward and backward amplitudes and phases; Fig. 8 shows the effects of the white-noise in the simulated reference response for the shaft's displacement inside the bearing 1.

The search domain was constrained in an interval from 25% to 200% of $Q_T$. The 3x3 grid of initial guesses was defined with values of $0.7 \times Q_T$, $Q_T$ and $1.3 \times Q_T$. The number of iterations is limited to 15 from each initial guess. The Jacobian used in the Newton-Raphson method (NRM) was approximated by central finite difference with a spacing of 2% of $Q_T$.

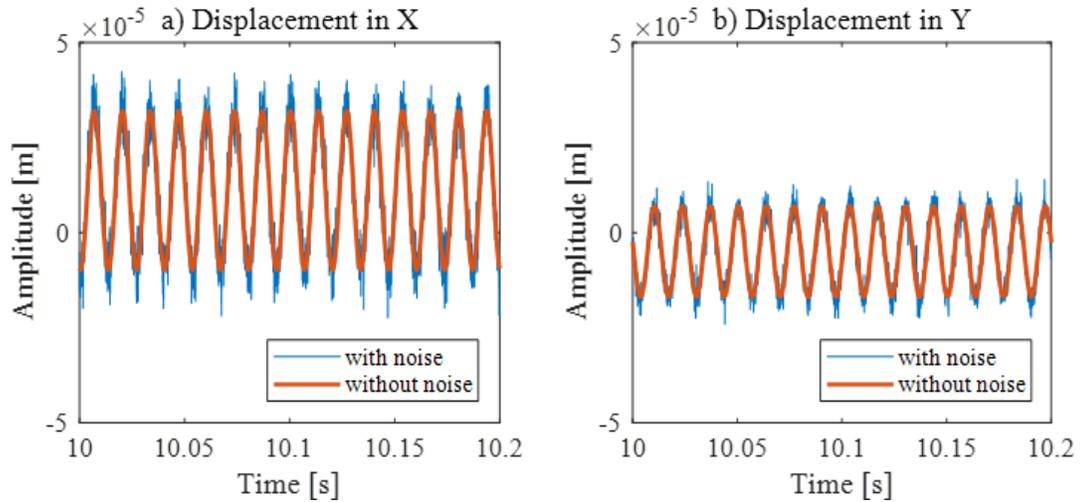

Fig. 8. Simulated displacements with and without noise at the center of bearing 1 with $Q_1 = Q_2 =$ 409.5 ml/min.

*3.2.1 Identification with rotor response evaluated at bearing's center*

This section presents the identification results using the shaft's displacement at the center of the bearings 1 and 2. Table 4 shows the relative error between identified and actual flowrate in bearing 1 ($Q_{1_{ident}}$ and $Q_{1_{ref}}$, respectively) for all the 64 cases, while the relative error between identified and actual flowrate in bearing 2 ($Q_{2_{ident}}$ and $Q_{2_{ref}}$, respectively) is shown in Table 5. Based on the obtained results, the proposed method allowed a satisfactory identification, with an average error of 0.52% and maximum error lower than 3%. In addition, small errors are observed over the entire flowrate range tested, which means that the method is applicable under flooded and starved lubrication conditions.

Table 4

Relative error of the identified oil flowrate in bearing 1, using vibration from the bearings center.

| Error $Q_{1_{ident}}$ (%) | | $Q_{2_{ref}}$ (ml/min) | | | | | | |
|---|---|---|---|---|---|---|---|---|
| | | 409.5 | 448.5 | 487.5 | 526.5 | 565.5 | 604.5 | 643.5 | 682.5 |
| $Q_{1_{ref}}$ (ml/min) | 409.5 | -0,27 | -0,04 | -0,81 | -0,27 | -0,97 | -1,99 | 0,05 | -0,67 |
| | 448.5 | 0,88 | 0,63 | 0,27 | -0,27 | 0,28 | -0,37 | 0,89 | 0,43 |
| | 487.5 | 0,10 | 0,01 | -0,13 | -0,36 | -0,46 | 0,06 | 0,78 | -0,12 |
| | 526.5 | -0,22 | -0,22 | -0,60 | 0,04 | -0,40 | -0,40 | -0,34 | -0,14 |
| | 565.5 | 0,04 | -0,34 | -0,16 | -0,10 | 0,39 | -0,15 | 0,07 | 0,28 |
| | 604.5 | 0,23 | -1,32 | 0,08 | 1,04 | 1,13 | -0,36 | 1,62 | -0,29 |
| | 643.5 | 1,29 | -1,39 | 0,58 | 1,04 | -0,78 | -0,48 | 0,87 | 0,98 |
| | 682.5 | 0,44 | -0,30 | -1,07 | -1,87 | -2,67 | 0,58 | -1,38 | 1,24 |

Table 5

Relative error of the identified oil flowrate in bearing 2, using vibration from the bearings center.

| Error $Q_{2ident}$ (%) | $Q_{2ref}$ (ml/min) | | | | | | | |
|---|---|---|---|---|---|---|---|---|
| $Q_{1ref}$ (ml/min) | 409.5 | 448.5 | 487.5 | 526.5 | 565.5 | 604.5 | 643.5 | 682.5 |
| 409.5 | -0,48 | -0,07 | -0,08 | -0,31 | -0,10 | 0,33 | 0,39 | -0,27 |
| 448.5 | 0,70 | 0,62 | 0,13 | 0,09 | -0,17 | 1,50 | 0,79 | 0,28 |
| 487.5 | -0,09 | 0,16 | -0,22 | 0,29 | 0,26 | 0,68 | -0,22 | -1,08 |
| 526.5 | 0,16 | -0,13 | 0,17 | 0,24 | 0,41 | -0,10 | 1,24 | 0,40 |
| 565.5 | -0,18 | -0,29 | 0,02 | -0,75 | 0,87 | 0,70 | 0,03 | -1,37 |
| 604.5 | -0,15 | 0,74 | -0,15 | -0,10 | 0,45 | -0,10 | 1,80 | -1,02 |
| 643.5 | 0,45 | 0,36 | -0,12 | -0,15 | -0,52 | -0,09 | 0,61 | -0,68 |
| 682.5 | -0,77 | 0,13 | 0,17 | -0,31 | -0,93 | 0,42 | 0,91 | 0,78 |

Table 6 shows the computation time, in hours, of the identification process for each case. Using a computer with an Intel Core i7-7700 3.6 GHz and 16 Gb of RAM, the average computation time is 1.2 h and the average number of iterations is 3. The costliest case ($Q_{1ref}$ = 448.5 ml/min and $Q_{2ref}$ = 604.5 ml/min) required 10 iterations, spending 3.9 h. The search path in this case is shown in Fig.9, where the blue and red lines show the roots curves of error functions in bearing 1 and bearing 2, respectively. From the Fig. 9, it is possible to note that the slower convergence occurred because the search oscillates around the root curves, thus, a damped method may improve the identification efficiency in this case. Despite that, the method needed just the first initial guess to identify the oil flowrates in all cases.

Table 6

Identification processing time using vibration from the bearings center.

| Processing Time (h) | $Q_{2ref}$ (ml/min) | | | | | | | |
|---|---|---|---|---|---|---|---|---|
| $Q_{1ref}$ (ml/min) | 409.5 | 448.5 | 487.5 | 526.5 | 565.5 | 604.5 | 643.5 | 682.5 |
| 409.5 | 1,2 | 1,3 | 1,4 | 1,4 | 2,1 | 1,4 | 1,0 | 1,0 |
| 448.5 | 1,6 | 1,2 | 2,1 | 2,4 | 1,7 | 3,9 | 1,4 | 1,3 |
| 487.5 | 1,7 | 1,6 | 0,9 | 1,1 | 1,5 | 2,4 | 2,4 | 1,3 |
| 526.5 | 1,3 | 1,8 | 0,9 | 0,5 | 0,9 | 1,3 | 1,0 | 1,0 |
| 565.5 | 2,4 | 1,9 | 1,6 | 0,8 | 0,5 | 0,8 | 0,9 | 1,0 |
| 604.5 | 1,8 | 1,7 | 2,9 | 1,0 | 0,9 | 1,2 | 1,2 | 1,5 |
| 643.5 | 1,4 | 1,7 | 2,4 | 1,0 | 0,9 | 1,9 | 0,8 | 0,8 |
| 682.5 | 1,0 | 1,7 | 2,7 | 1,0 | 0,9 | 1,9 | 0,8 | 0,5 |

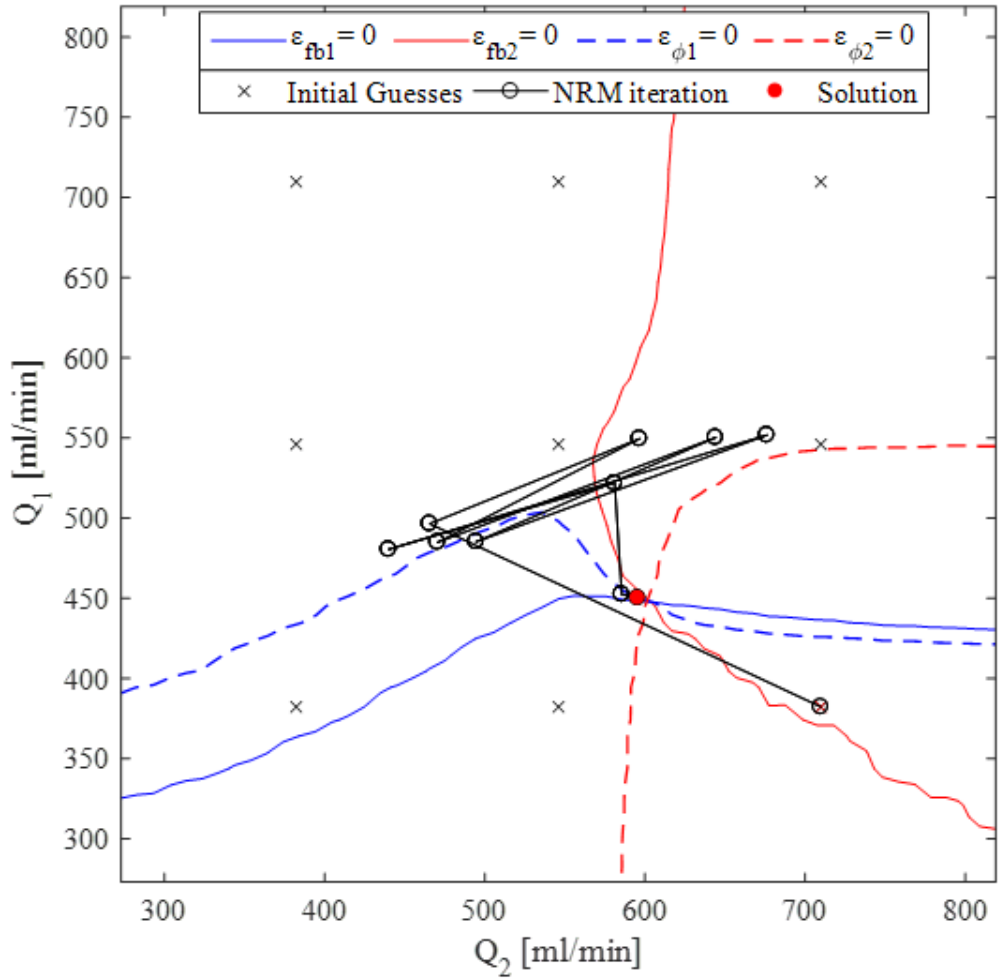

Fig. 9. Identification search in the case of $Q_{1ref} = 448.5$ ml/min and $Q_{2ref} = 604.5$ ml/min.

*3.2.2 Identification with rotor response evaluated outside the bearings*

In practical situations, the shaft displacement cannot always be measure at the bearing center (inside the bearing), hence the proposed identification method was also tested using the rotor response measured outside the bearings. Tables 7 and 8 show the relative error in the oil flowrates identified in bearings 1 and 2, respectively, using the responses at nodes 5 and 21, coincident with the bearings' side towards the shaft extremity.

Table 7

Relative error of the identified oil flowrate in bearing 1, using vibration from outside the bearings.

| Error $Q_{1ident}$ (%) $Q_{1ref}$ (ml/min) | $Q_{2ref}$ (ml/min) | | | | | | | |
|---|---|---|---|---|---|---|---|---|
| | 409.5 | 448.5 | 487.5 | 526.5 | 565.5 | 604.5 | 643.5 | 682.5 |
| 409.5 | -0,36 | 0,06 | -0,19 | -0,66 | -0,18 | -0,32 | -0,26 | -0,07 |
| 448.5 | 0,77 | 0,31 | -0,12 | 0,08 | -0,24 | 0,33 | 0,08 | -0,15 |
| 487.5 | -0,10 | -0,04 | 0,03 | -0,01 | -0,10 | -0,06 | -0,32 | -0,06 |
| 526.5 | -0,22 | 0,08 | 0,21 | -0,16 | -0,25 | 0,23 | 0,91 | 0,01 |
| 565.5 | -0,17 | -0,70 | -0,60 | 0,51 | 0,04 | 0,46 | -0,09 | 0,59 |
| 604.5 | -0,14 | -0,15 | -0,06 | -0,06 | 0,41 | -0,65 | 0,21 | 0,76 |
| 643.5 | 0,00 | 0,46 | 0,23 | 0,11 | 1,20 | 0,03 | 1,05 | 1,23 |
| 682.5 | -0,54 | -0,98 | -0,61 | -0,48 | -0,57 | 0,39 | -0,20 | -1,25 |

Table 8

Relative error of the identified oil flowrate in bearing 2, using vibration from outside the bearings.

| Error $Q_{2ident}$ (%) $Q_{1ref}$ (ml/min) | $Q_{2ref}$ (ml/min) | | | | | | | |
|---|---|---|---|---|---|---|---|---|
| | 409.5 | 448.5 | 487.5 | 526.5 | 565.5 | 604.5 | 643.5 | 682.5 |
| 409.5 | -0,44 | 0,25 | -0,07 | -0,21 | -0,11 | 0,03 | -0,36 | -0,70 |
| 448.5 | -0,02 | 0,29 | 0,01 | -0,24 | 0,32 | -0,48 | -1,63 | -1,29 |
| 487.5 | -0,10 | 0,30 | -0,02 | 0,09 | 0,24 | 0,03 | 1,08 | -0,56 |
| 526.5 | -0,57 | -0,45 | 0,02 | 0,68 | -0,26 | 0,98 | -0,10 | 0,94 |
| 565.5 | -0,55 | 0,80 | 0,10 | 0,00 | 0,06 | -0,01 | -0,13 | -1,97 |
| 604.5 | -0,25 | 0,37 | -0,06 | 0,04 | 0,33 | 0,58 | -0,52 | -1,24 |
| 643.5 | -0,08 | 0,00 | -0,04 | -0,04 | 0,43 | -0,02 | 0,61 | -1,11 |
| 682.5 | -0,36 | 0,14 | -0,04 | -0,14 | -0,01 | 0,14 | 0,56 | 0,19 |

According to the results obtained, the oil flowrates are successfully identified with an average error of 0.36% and a maximum error lower than 2%. However, the average number of iterations increased to 10, reaching 70 iterations in some cases, which means that all the initial guesses were tried before finding the root. The number of iterations increases due to the presence of multiple roots of some error functions used in the identification, mainly $\varepsilon_{\varphi 1}$ and $\varepsilon_{\varphi 2}$ (Fig. 10). In addition, the search diverges from some initial guesses (overshoots the oil flowrate bounds) and a new root search is then restarted from another initial guess, also increasing the number of iterations (Fig 11).

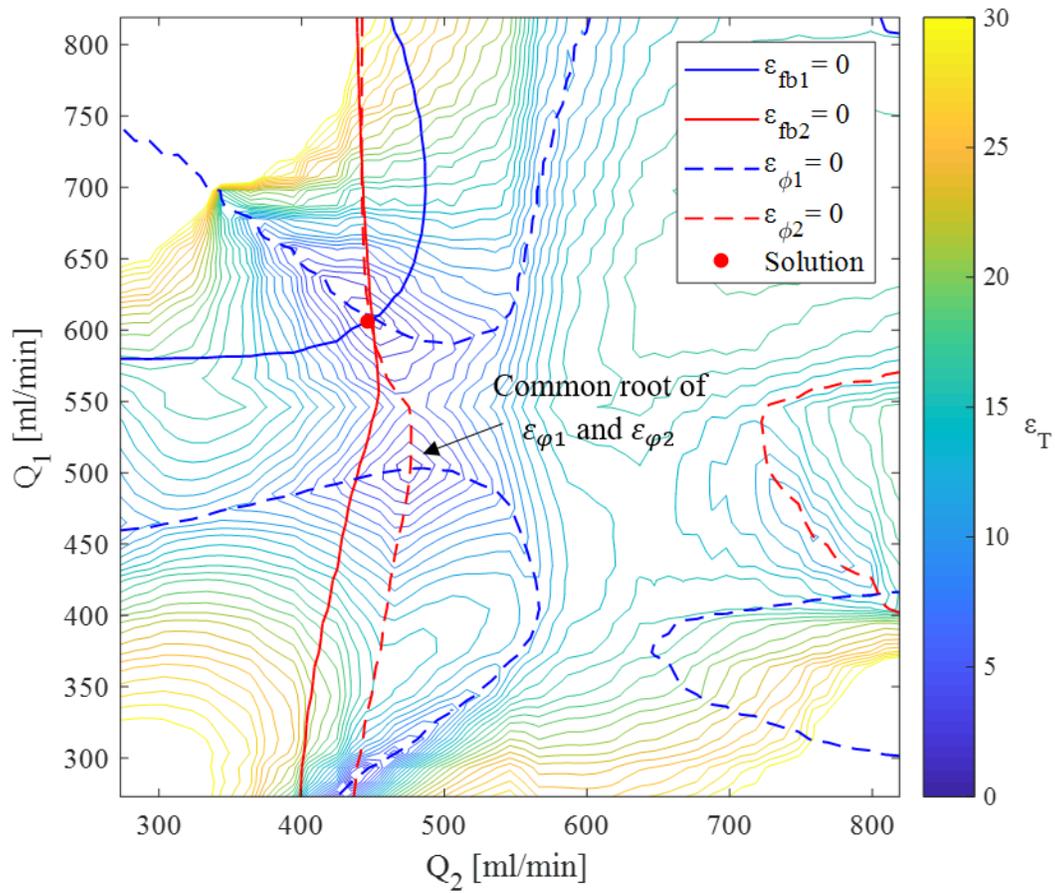

Fig. 10. Contour of the total error ($\varepsilon_T$) and root curves of error functions in a costly case ($Q_{1ref}$ = 604.5 ml/min and $Q_{2ref}$ = 488.5 ml/min).

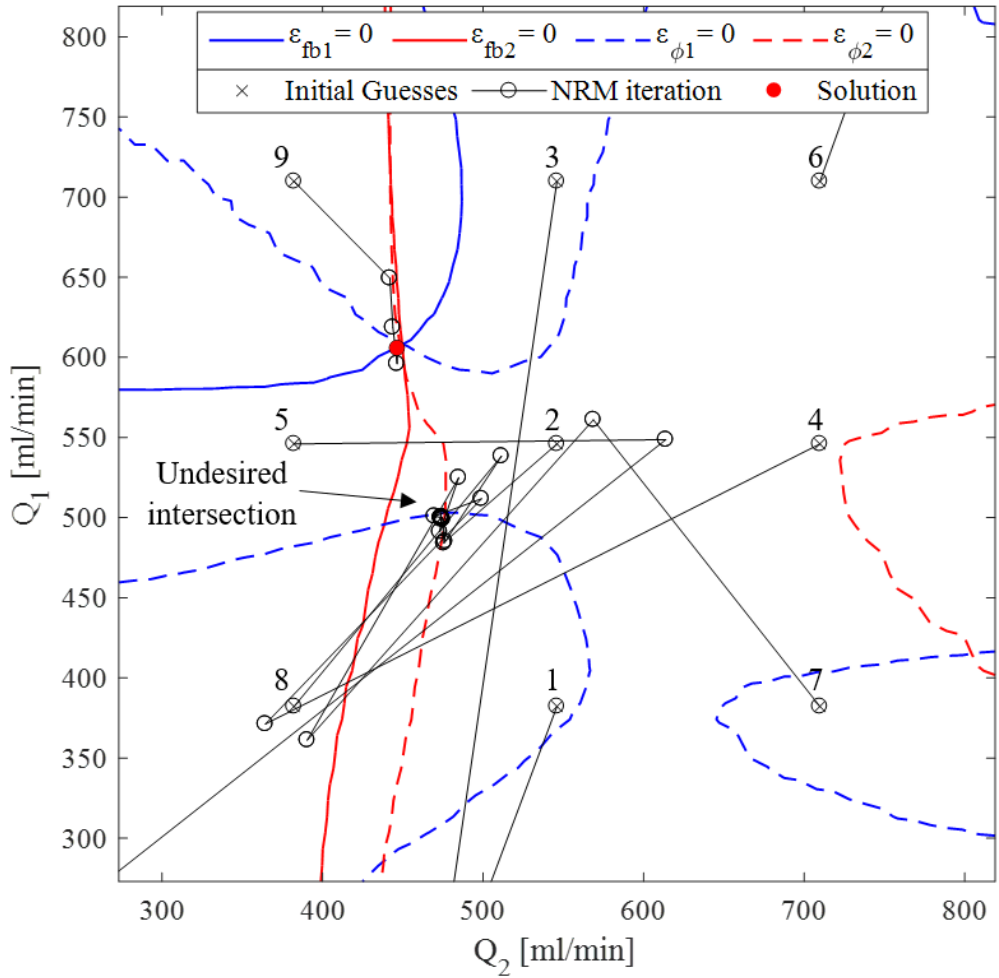

Fig. 11. Identification search in a costly case ($Q_{1ref}$ = 604.5 ml/min and $Q_{2ref}$ = 488.5 ml/min). Initial guesses numbered in order of attempt.

Observing the Fig. 10, it is possible to visualize the intersection of all error functions roots (solution), but there is also an undesired intersection between the $\varepsilon_{\varphi 1}$ and $\varepsilon_{\varphi 2}$ roots that results a higher total error. As presented in Fig. 11, the search tracks the intersection of $\varphi$ error roots, but, at this point, the total error ($\varepsilon_T$) is not low enough to satisfy the stop criterion, forcing a new attempt with another initial guess. Observing the tried initial guesses (Fig. 11), it is possible to note that all guesses were needed because the search overshoots the oil flowrate bounds when started from the initial guesses 1, 3, 5 and 6, and got stuck in the wrong intersection when it started from initial guesses 2, 4, 7 and 8.

## 4. Conclusion

This paper presents a fault identification and diagnosis method based on the rotor vibration signals to identify the lubrication condition of hydrodynamic bearings. In view of this, a new approach was proposed to include the oil supply flowrate as an input parameter of a mass-conserving lubrication model. Then, a model-based identification process was developed to iteratively estimate the oil flowrate in the bearings of a rotating system.

Computational simulations showed that forward and backward amplitudes, as well as the orbit angle of the rotor vibration are sensitive to oil flowrate variation. From the results obtained, it is concluded that using the ratio between forward and backward amplitudes is more suitable than using them separately in the proposed identification method.

Numerical identification tests showed that the proposed method is capable to satisfactorily identify bearing oil supply flowrate under starved and flooded lubrication conditions. Even from highly noisy vibration data, the identified flowrates presented errors always below 3%. Moreover, the position of the data acquisition barely affected the errors, but showed a notable influence on the method efficiency. In some cases, evaluating the vibration signals measured outside the bearing, the number of iterations required for the identification increases drastically in comparison with the identification performed with the vibrations signals measured at the center of the bearings.

The novel method proposed to identify bearings lubrication condition has practical application as a monitoring technique for industrial machines, representing an important contribution to the field of machinery health monitoring. It is important to highlight that the identification process of the oil supply flowrates has been successfully performed even on noisy signals measured outside the bearing, indicating a promising tool to monitor rotating machines operating in real applications.

*CRediT authorship contribution statement*

M.V.M Oliveira: Conceptualization, Methodology, Software, Validation, Formal analysis, Investigation, Writing - original draft, Writing - review & editing, Visualization. B.Z. Cunha: Methodology, Software, Validation, Formal analysis, Investigation, Writing - review & editing. G.B. Daniel: Methodology, Validation, Formal analysis, Investigation, Writing - review & editing, Supervision, Project administration, Funding acquisition.

*Acknowledgments*

The authors would like to acknowledge CNPq and CAPES for the financial support. Special thanks to Matheus Victor Inacio for his insights into the root search problem.

*References*

[1] Chen J, Patton RJ. Robust model-based fault diagnosis for dynamic systems. Massachusetts, USA: Kluwer Academic Publishers; 1999.

[2] Bachschmid N, Pennachi P, Vania A. Identification of multiple faults in rotor systems. Journal of Sound and Vibration 2002; 254(2):327-366

[3] Alves DS, Daniel GB, Castro HF, Machado TH, Cavalca KT, Gecgel O, Dias JP, Ekwaro-Osire S. Uncertainty quantification in deep convolutional neural network diagnostics of journal bearings with ovalization fault. Mechanism and Machine Theory 2020; 149:103835.

[4] Papadopoulos CA, Nikolakopoulos PG, Gounaris GD. Identification of clearances and stability analysis for a rotor-journal bearing system. Mechanism and Machine Theory 2008; 43:411-426.

[5] Gertzos KP, Nikolakopoulos PG, Chasalevris AC, Papadopoulos CA. Wear identification in rotor-bearing systems by measurements of dynamic bearing characteristics. Computers and Structures 2011; 89:55–66.

[6] Machado TH, Cavalca KL. Investigation on an experimental approach to evaluate a wear model for hydrodynamic cylindrical bearings systems. Applied Mathematical Modeling 2016; 40:9546–9564.

[7] Mendes RU, Machado TH, Cavalca KL. Experimental wear parameters identification in hydrodynamic bearings via model based methodology. Wear 2017; 372–373:116–29.

[8] Alves DS, Fieux G, Machado TH, Keogh PS, Cavalca KL. A parametric model to identify hydrodynamic bearing wear at a single rotating speed. Tribology International 2021; 153:106640.

[9] Bonneau D, Frene J. Film formation and flow characteristics at the inlet of a starved contact—theoretical study. Journal of Lubrication Technology 1983; 105(2):178-185.

[10] Artiles A, Hesmat H. Analysis of starved journal bearings including temperature and cavitation effects. Journal of tribology 1985; 107(1):1-13.

[11] Vincent B, Maspeyrot P, Frene J. Starvation and cavitation effects in finite grooved journal bearing. Tribology Series 1995; 30:455-464.

[12] Tanaka M. Journal bearing performance under starved lubrication. Tribology International 2000; 33(3-4):259-264.

[13] Vijayaraghavan D, Keith Jr TG, Brewe DE. Effect of lubricant supply starvation on the thermohydrodynamic performance of a journal bearing. Tribology Transactions 1996; 39(3):645-653.

[14] Hashimoto H, Ochiai M. Stabilization Method for Small-Bore Journal Bearings Utilizing Starved Lubrication. Journal of Tribology 2010; 132(1):011703.

[15] Guo F, Zang S, Li C, Wong PL, Guo L. Lubrication film generation in slider-on-disc contact under limited lubricant supply. Tribology International 2018; 125:200-208.

[16] Liu CL, Guo F, Wong PL. Characterisation of starved hydrodynamic lubricating films. Tribology International 2019; 131:694-701.

[17] Poddar S, Tandon N. Study of Oil Starvation in Journal Bearing Using Acoustic Emission and Vibration Measurement Techniques. Journal of Tribology 2020; 142(12):121801.

[18] Poddar S, Tandon N. Classification and detection of cavitation, particle contamination and oil starvation in journal bearing through machine learning approach using acoustic emission signals. Proceedings of the Institution of Mechanical Engineers, Part J: Journal of Engineering Tribology 2021.

[19] Ausas RF, Jai M, Buscaglia GC. A mass-conserving algorithm for dynamical lubrication problems with cavitation. Journal of Tribology 2009; 131(3):031702.

[20] Lund JW. Spring and damping coefficients for the tilting-pad journal bearing. ASLE Transactions 1964; 7(4):342-352.

[21] Nelson HD, McVaugh JM. The Dynamics of Rotor-Bearing Systems Using Finite Elements. Journal of Engineering for Industry 1976; 98(2):593-600.

[22] Lee CW, Han YS. The Directional Wigner Distribution and Its Applications. Journal of Sound and Vibration 1998; 216(4):585–600.